

\baselineskip=14pt
\parskip=10pt
\def\halmos{\hbox{\vrule height0.15cm width0.01cm\vbox{\hrule height
  0.01cm width0.2cm \vskip0.15cm \hrule height 0.01cm width0.2cm}\vrule
  height0.15cm width 0.01cm}}
\font\eightrm=cmr8 
\font\eighttt=cmtt8
\magnification=\magstephalf

\def\1{{\overline{1}}}
\def\2{{\overline{2}}}
\parindent=0pt
\overfullrule=0in

\def\frac#1#2{{#1 \over #2}}
\bf
\centerline
{
Generalizing and Implementing  
}
\centerline
{
Michael Hirschhorn's AMAZING Algorithm for Proving Ramanujan-Type Congruences\footnote{$^1$}
{
\eightrm
June 27, 2013. Exclusively published in the Personal Journal of Shalosh B. Ekhad and
Doron Zeilberger as well as in arxiv.org . Accompanied by Maple packages
{\eighttt HIRSCHHORN} and {\eighttt BOYLAN} downloadable from the front of this article
\hfill\break
{\eighttt http://www.math.rutgers.edu/\~{}zeilberg/mamarim/mamarimhtml/mh.html},  where
there are several sample input and output files containing computer-generated theorems (and proofs!). \hfill \break
The work of DZ was supported in part by the National Science Foundation of the United States of America.
}
}
\rm
\bigskip
\centerline{ Edinah GNANG and Doron  ZEILBERGER}

{\bf Preamble} 

Let $p(n)$ be the number of integer partitions of $n$. Euler famously proved that
$$
\sum_{n=0}^{\infty} p(n) q^n = \prod_{i=1}^{\infty} \frac{1}{1-q^i}  \quad .
$$
Srinivasa Ramanujan famously discovered (by glancing at a table of $p(n)$ for $1 \leq n \leq 200$,
computed by the {\it analytic machine}, Major Percy Alexander MacMahon's head) the three congruences
$$
p(5m+4)   \equiv 0  \pmod  5  \quad ,
$$
$$
p(7m+5)  \equiv 0  \pmod  7  \quad ,
$$
$$
p(11m+6)  \equiv 0  \pmod  {11}  \quad .
$$

The first two are really easy, and the proofs that G.H. Hardy chose to present in his
classic book ``Ramanujan'' ([Ha], pp. 87-88), slightly streamlined,
go as follows.

First recall the (purely elementary and {\it shaloshable}) identities of Euler and Jacobi :
$$
E(q)=\prod_{i=1}^{\infty} (1-q^i) =\sum_{n=-\infty}^{\infty} (-1)^n q^{(3n^2+n)/2} \quad , \quad and
$$
$$
E(q)^3=\sum_{n=0}^{\infty} (-1)^n (2n+1) q^{(n^2+n)/2} \quad .
$$

Also recall the obvious fact (but {\it extremely useful} [e.g. the AKS algorithm!] ),
that follows  from the binomial theorem and Fermat's little theorem, that for every prime $\ell$, and any polynomial,
or formal power series, $f(q)$, $f(q)^\ell \equiv f(q^\ell) \, \pmod  \ell$ . In particular
$E(q)^\ell \equiv E(q^\ell) \, \pmod  \ell$ .

{\bf p(5n+4) is divisible by 5}

Since $\{(n^2+n)/2 \,\, mod \,\, 5 \, \, ; \, \, 0 \leq n \leq 4 \,\, , \,\, 2n+1 \not \equiv 0 \pmod 5 \}=\{0,1\}$, we have:
$$
E(q)^3  \equiv J_0+ J_1 \pmod 5 \quad,
$$
where $J_i$ consists of those terms in which the power of $q$ is congruent to $i$ modulo $5$.
Now
$$
\sum_{n=0}^{\infty} p(n) q^n =
E(q)^{-1} =\frac{(E(q)^3)^3}{E(q)^{10}} =\frac{(E(q)^3)^3}{(E(q)^5)^2} \equiv \frac{(J_0+J_1)^3}{ E(q^5)^2}  \pmod 5 \quad .
$$
Since $(J_0+J_1)^3=J_0^3+3J_0^2J_1 + 3J_0J_1^2 +J_1^3$, whose terms consist of powers 
of $q$ that are $0,1,2,3$ modulo $5$, respectively, none of the powers of $q$ that are congruent to
$4$ modulo $5$ show up, and  hence the coefficient of $q^{5n+4}$ is always $0$ modulo $5$. \halmos

{\bf p(7n+5) is divisible by 7}

Since $\{(n^2+n)/2 \,\, mod \,\, 7 \, \, ; \, \, 0 \leq n \leq 6 \,\, , \,\, 2n+1 \not \equiv 0 \pmod 7 \}=\{0, 1, 3\}$, we have:
$$
E(q)^3 \equiv J_0+ J_1 +J_3 \pmod 7 \quad,
$$
where $J_i$ consists of those terms in which the power of $q$ is congruent to $i$ modulo $7$.
Now
$$
\sum_{n=0}^{\infty} p(n) q^n =
E(q)^{-1} = \frac{(E(q)^3)^2}{E(q)^7} \equiv \frac{(J_0+J_1+J_3)^2}{E(q^7)}  \pmod 7 \quad ,
$$
Since $(J_0+J_1+J_3)^2=J_0^2+J_1^2 + J_3^2 + 2J_0J_1+2J_0J_3+2J_1J_3$, whose terms consist of powers 
of $q$ that are $0,2,6,1,3,4$ modulo $7$, respectively, none of the powers of $q$ that are congruent to
$5$ modulo $7$ show up, and  hence the coefficient of $q^{7n+5}$ is always $0$ modulo $7$. \halmos

At the bottom of page 88 of Hardy's above-mentioned classic ``Ramanujan''[Ha], he states

``{\it There does not seem to be an equally simple proof that $p(11n+6)$ is divisible by $11$}''.

Over the years there were many  proofs,  but {\bf none} as
{\bf simple} and {\bf elementary} and, most importantly, {\bf beautiful!}, as the one recently found by Michael Hirschhorn [Hi].

{\bf Michael Hirschhorn's proof that  p(11n+6)  is divisible by 11}

The proof in [Hi] goes like this. It starts the same way.

Since $\{(n^2+n)/2 \,\, mod \,\, 11 \, \, ; \, \, 0 \leq n \leq 10 \,\, , \,\, 2n+1 \not \equiv 0 \pmod {11} \}=\{0, 1, 3, 6, 10\}$, we have:
$$
E(q)^3 \equiv J_0+ J_1 + J_3 + J_6 + J_{10}  \pmod {11} \quad,
$$
where $J_i$ consists of those terms in which the power of $q$ is congruent to $i$ modulo $11$.
Now
$$
\sum_{n=0}^{\infty} p(n) q^n =
E(q)^{-1} =\frac{(E(q)^3)^7}{E(q)^{22}} \equiv \frac{(J_0+ J_1 + J_3 + J_6 + J_{10} )^7}{E(q^{11})^2}  \pmod {11} \quad .
$$
Alas, now the part consisting of the powers that are congruent to $6$ modulo $11$ in the polynomial
$(J_0+ J_1 + J_3 + J_6 + J_{10} )^7 \pmod {11}$
is {\bf not} identically zero  modulo $11$, but a certain
polynomial of degree $7$ in  $\{J_0, J_1,  J_3,  J_6,  J_{10} \}$, (over $GF(11)$) let's call it $POL$.

It is readily seen that, introducing an auxiliary variable $t$, that
$$
POL(J_0, J_1,  J_3,  J_6,  J_{10})=
Coeff_{t^6} \left ( J_0+ J_1 t + J_3 t^3 + J_6 t^6+ J_{10} t^{10} \right )^7 \pmod {11} \pmod {t^{11}-1} \quad ,
$$
that is {\bf not} identically zero.

But, Since $\{(3n^2+n)/2 \,\, mod \,\, 11 \, \, ; \, \, 0 \leq n \leq 10 \, \}=\{0, 1, 2, 4, 5, 7\}$, we have,
$$
E(q)= E_0 + E_1 + E_2 +E_4 + E_5 + E_7 \quad  ,
$$
where $E_i$ consists of those terms in which the power of $q$ is congruent to $i$ modulo $11$,
and 
$$
(E(q)^3)^4=E(q)^{12}= E(q)^{11} E(q) \equiv  E(q^{11}) E(q) \pmod {11} \quad,
$$
so
$$
(J_0+J_1+J_3+J_6+J_{10})^4  \equiv  E(q^{11})   (E_0 + E_1 + E_2 +E_4 + E_5 + E_7 ) \pmod {11} \quad .
$$
By expanding the left side and extracting the complementary powers $\pmod {11}$ ( $\{3, 6, 8, 9 ,10\}$),
we get five polynomials of degree $4$, let's call them $Q_3, Q_6, Q_8, Q_9, Q_{10}$
that we know are $0$ modulo $11$ (once the $J_i$'s are replaced by the formal power series they stand for).
For $m \in \{3,6,8,9,10\}$, we have
$$
Q_{m}(J_0, J_1,  J_3,  J_6,  J_{10})=
Coeff_{t^m} \left ( J_0+ J_1 t + J_3 t^3 + J_6 t^6+ J_{10} t^{10} \right )^4 \pmod {11} \pmod {t^{11}-1} \quad .
$$
Then we ask our beloved computer to find five polynomials 
of degree $3$, (in the variables  $\{J_0, J_1,  J_3,  J_6,  J_{10} \}$),
let's call them $R_3, R_6, R_8, R_9, R_{10}$, such that
$$
POL \equiv  R_3 Q_3 + R_6 Q_6 + R_8 Q_8 + R_9 Q_9 + R_{10} Q_{10} \pmod {11} \quad .
$$
Since it succeeded ({\it a priori} there was no guarantee!), we are done!!
{\it Quod Erat Demonstratum}. \halmos

See the output file {\tt http://www.math.rutgers.edu/\~{}zeilberg/tokhniot/oHIRSCHHORN1v}, that
contains the above three proofs, (and four other ones!), that was generated, 
by running the Maple package {\tt HIRSCHHORN} (that accompanies this article), in three seconds!

{\bf More Ramanujan Type Congruences}

Let's consider, more generally,
$$
\sum_{n=0}^{\infty} p_{-a}(n) q^n = \prod_{i=1}^{\infty} \frac{1}{(1-q^i)^a}  \quad .
$$

(Note that $p_{-1}(n)=p(n)$ and $p_{24}(n)=\tau(n-1)$, where $\tau(n)$ is Ramanujan's $\tau$-function).

There are many known Ramanujan-type congruences for $p_{-a}(n)$. Matthew Boylan [B] (Theorem 1.3,
where our $p_{-a}(n)$ is denoted by $p_a(n)$, and the entry $r=27, l=31$ is erroneous)
has found all of them for $a$ odd and $\leq 47$. 

The first few are (here we restricted our search to primes $\geq 2a+1$).
$$
p_{-1}(5n+4) \equiv 0  \pmod  {5}  \quad , \quad 
p_{-1}(7n+5) \equiv 0  \pmod  {7}  \quad , \quad 
p_{-1}(11n+6) \equiv 0  \pmod  {11}  \quad , \quad 
$$
(Ramanujan's)
$$
p_{-2}(5n+2) \equiv 0  \pmod  {5}  \quad , \quad 
p_{-2}(5n+3) \equiv 0  \pmod  {5}  \quad , \quad 
p_{-2}(5n+4) \equiv 0  \pmod  {5}  \quad , \quad 
$$
$$
p_{-3}(11n+7) \equiv 0  \pmod  {11}  \quad , \quad 
p_{-3}(17n+15) \equiv 0  \pmod  {17}  \quad , \quad 
$$
$$
p_{-5}(11n+8) \equiv 0  \pmod  {11}  \quad , \quad 
p_{-5}(23n+5) \equiv 0  \pmod  {23}  \quad , \quad 
$$
$$
p_{-7}(19n+9) \equiv 0  \pmod  {19}  \quad , \quad 
$$
$$
p_{-9}(19n+17) \equiv 0  \pmod  {19}  \quad , \quad 
p_{-9}(23n+9) \equiv 0  \pmod  {23}  \quad , \quad 
$$
$$
p_{-21}(47n+42) \equiv 0  \pmod  {47}  \quad .
$$

Thanks to the impressive algorithm of Silviu Radu[R1], every such congruence (and even more general ones, see [R1]),
is effectively (and fairly efficiently!) decidable. Let's hope that Radu would post a public implementation
of his method.  Since no such an implementation seems to exist, we Emailed Radu, who
kindly[R2] showed us how to deduce these (except for the last two, that we are sure
can be done just as easily) from his powerful algorithm, by specifying the $N_0$ for which
checking them for $0 \leq n \leq N_0$ would imply them for {\it all}  $\,\, 0 \leq n < \infty$.

As impressive as Radu's algorithm is, it is {\bf not} elementary. It uses the `fancy', and intimidating,
theory of modular functions, that being analytic, is not
quite legitimate according to our finitistic and discrete philosophy of mathematics. Hence it is still interesting (at least to us!)
to find elementary, Hirschhorn-style proofs. Also, by the {\it principle of serendipity} our
extension and implementation of Hirschhorn's method may lead to new things that even
modular functions can {\bf not} do.

{\bf Extending Hirschhorn's Method}

Suppose that, for some prime $\ell$ and some integer $r$ ($0 \leq r<\ell$),  we want to prove a congruence of the type
$$
p_{-a}(\ell n+r) \equiv 0  \pmod  {\ell}   \quad .
$$

We first find the smallest integer $\alpha$ such that $b:=(\alpha \ell-a)/3$ is an integer, noting that
$$
E(q)^{-a} =\frac{(E(q)^3)^b}{E(q)^{\alpha \ell}} \equiv
 \frac{(E(q)^3)^b}{E(q^\ell)^\alpha} \pmod {\ell} \quad .
$$
We now define the subset of $\{0,1, \dots, \ell-1\}$:
$$
Jset(\ell):=
\{(n^2+n)/2 \,\, mod \,\, \ell \, \, ; \, \, 0 \leq n \leq \ell -1 \,\, , \,\, 2n+1 \not \equiv 0 \pmod \ell \} \quad ,
$$
and write
$$
E(q)^3 \equiv \sum_{i \in  Jset(\ell)} J_i   \pmod {\ell} \quad,
$$
where $J_i$ consists of those terms in which the power of $q$ is congruent to $i$ modulo $\ell$.
Next we define $POL$ to be the  polynomial, in the set of variables $\{ J_i \, ; \, i \in Jset(\ell) \}$,
$$
POL(\{ J_i \, ; \,  i \in Jset(\ell) \})=
Coeff_{t^r} \left [  \, \left ( \, \sum_{i \in  Jset(\ell)} \, J_i t^i \right )^b \right] \pmod {\ell} \pmod {t^\ell-1} \quad .
$$

Now, if we are {\bf lucky}, the polynomial $POL(\{ J_i \})$ would be identically zero (modulo $\ell$). In that case
we have a {\it Ramanujan-style} proof, since the powers of $q$ that are congruent to $r$ modulo $\ell$
in $((E(q)^3)^{b}$, and hence in $E(q)^{-a}$, do not show up!

Otherwise, we need to resort to {\it Hirschhorn's} enhancement.

Analogously to $Jset(p)$, let's define 
$$
Eset(\ell):=\{ (3n^2+n)/2 \,\, mod \,\, \ell \, \, ; \,\, 0 \leq n \leq \ell -1 \} \quad,
$$
the set of residue classes modulo $\ell$ that show up as powers in the sparse Euler Pentagonal Theorem
expression for $E(q)$.

Now let $c$ be the reciprocal of $3$ modulo $\ell$, and let $d=(3c-1)/\ell$. Then
$$
(E(q)^3)^c=E(q)E(q)^{3c-1}=E(q)E(q)^{d\ell} \equiv E(q) (E(q^\ell))^d \pmod {\ell} .
$$
Now define a set of polynomials, for each $0 \leq m<\ell$ that is {\bf not} in $Eset(\ell)$
(i.e. for the members of the complement of $Eset(\ell)$):
$$
Q_m:=
Coeff_{t^m} \left [ \, \left ( \, \sum_{i \in  Jset(\ell)} J_i t^i \, \right)^c \right ] \pmod {\ell} \pmod {t^\ell - 1} \quad , \quad
m \not\in Eset(\ell) \quad .
$$
We know that all the $Q_m(\{J_i\})$ [$m \not \in Eset(\ell)$] are $0$ modulo $\ell$ (once the $J_i's$ are replaced by
the formal power series, in $q$, that they stand for).

Finally, we decide whether the polynomial $POL$ (that lives in the polynomial ring over the Galois 
Field $GF(\ell)$ in the $J_i$'s), or one of its powers, belongs to the {\bf ideal} generated by the polynomials
$Q_m$. This can be done (for small $\ell$) either directly, using
{\it undetermined coefficients}, and for larger $\ell$, using the {\it Buchberger algorithm} (alias {\it Gr\"obner bases}).

{\bf The Big Disappointment}

We naively hoped that Hirschhorn's method, as explicated and generalized above, would work for all of these
other congruences. To our dismay, it failed to prove the congruence $p_{-3}(17n+15) \equiv 0 \pmod {17}$.

It turns our that for the specialization 
$$
J_0=1 \, , \, J_1=1  \, , \, J_3 =2 \, , \, J_4=10 \, , \, J_6=9 \, , \, J_{10}=11 \, ,\, J_{11}=15 \, , \, J_{15}=12 \quad ,
$$
all the $Q_m$ are zero (modulo $17$) {\bf but} $POL \equiv 6 \pmod {17} \neq 0$. So, of course, $POL$ is not
in the ideal generated by the $Q_m$ in $GF(17)[J_0,J_1,J_3,J_4,J_6,J_{10},J_{11},J_{15}]$.

{\bf But there is Hope}

The Euler and Jacobi identities are but the first two in an {\it infinite} sequence of identities,
the  {\it Macdonald identities}[M] made famous  in Freeman Dyson's[D] historic 1972 Gibbs Lecture.

In fact, the next-in-line in Macdonald's identities, earlier found by Winquist[W], was already
used to give ``a Ramanujan-style proof'' of $p(11m+6) \equiv 0 \pmod {11}$. We strongly believe
that {\it every Ramanujan-type congruence} that can be proved using Radu's[R1]  beautiful algorithm
(that relies on the theory of modular functions), has either a ``Ramanujan-style'', or ``Hirschhorn-style''
proof, by using one of the Macdonald identities, that in spite of their ``fancy'' pedigree
(Lie theory) are {\bf purely elementary}.

{\bf The Maple package HIRSCHHORN}

Everything (and more) is implemented in the Maple package {\tt HIRSCHHORN} available from

{\tt http://www.math.rutgers.edu/\~{}zeilberg/tokhniot/HIRSCHHORN} \quad .

The webpage

{\tt http://www.math.rutgers.edu/\~{}zeilberg/mamarmim/mamarimhtml/mh.html}

contains several computer-generated articles outputted by that package.

{\bf Gr\"obner via special cases}

For $\ell >11$, both $POL$ and the $\{ Q_m \}$ get too big for Maple. But by
doing sufficiently many specializations (mod $\ell$) for 
a subset of the variables $J_i's$ one
can get a {\bf fully rigorous} proof of ideal membership. See 
procedures {\tt TerseMikeProof,  TerseMikeProofG, TerseMikeProofGviaSC},
that use, respectively, undetermined coefficients, Gro\"obner bases,
and  Gro\"obner bases via special cases. As we already pointed out above, we
are not always guaranteed success.

{\bf Future Directions}

We believe that our extension of Hirschhorn's method could be generalized to
more general $q$-series, including those that are {\it not} modular functions.

{\bf FIRST Encore: The Maple package BOYLAN}

The Maple package {\tt BOYLAN} available from

{\tt http://www.math.rutgers.edu/\~{}zeilberg/tokhniot/BOYLAN} 

reproduces and extends Theorem  1.3 of [Bo], albeit empirically.

See

{\tt http://www.math.rutgers.edu/\~{}zeilberg/tokhniot/oBOYLAN1} 

for a reproduction of the original (in less than two seconds), and

{\tt http://www.math.rutgers.edu/\~{}zeilberg/tokhniot/oBOYLAN2} 

for many more congruences (going as far as $a=399$).

{\bf SECOND Encore: Infinitely Many Congruences} ({\it All having Ramanujan-style proofs!})

Now that, thanks to Radu[R1], any {\it specific} congruence of the form $p_{-a}(\ell n+r) \equiv 0 \pmod {\ell}$,
is {\it purely routine} (or, more politely, {\it algorithmically provable}, or {\it shaloshable}), the
next stage would be to come up with ``infinitely many congruences''.

There is, of course, a {\it cheap} way to get ``infinitely many'' such congruences, namely when
$a=\ell-3$, since
$$
\frac{1}{E(q)^{\ell-3}} \equiv \frac{E(q)^3}{E(q^\ell)} \pmod \ell ,
$$
and since the set $Jset(\ell)$ is about  one half of all residue classes, we get many $r$'s 
(all the members of the complement of $Jset(\ell)$).

But, a little less trivially, we can generalize the Ramanujan proof of  $p(7n+5) \equiv 0 \pmod 7$,
to the following proposition (we hesitate to call it a theorem, for two reasons. First it is a bit
shallow, and second we only have a sketch of a proof [that we are sure can be easily completed, but
we have better things to do]).

{\bf Proposition}: Let $\ell$ be a prime that is either $7$ or $11$ modulo $12$ and let
$r:= (\ell-6)/24 \pmod {\ell}$, then
$$
p_{-(\ell-6)} (n \ell +r) \equiv 0 \pmod \ell \quad .
$$

{\bf Sketch of a  Ramanujan style proof}. It is easy to see that  $r:= (\ell-6)/24 \pmod {\ell} \not \in Jset(\ell)+Jset(\ell) \pmod \ell$,
thanks to the following (presumably elementary) lemma, that we verified empirically for $\ell \leq 2000$.

{\bf Elementary Lemma:}
Let $\ell$ be a prime that leaves remainder $7$  or $11$ when divided by $12$. Then for any $0 \leq n_1, n_2 <\ell$ such that
$$
\frac{n_1(n_1+1)}{2}+\frac{n_2(n_2+1)}{2} \equiv r \pmod {\ell} \quad,
$$
we must have either $n_1=(\ell -1)/2$ or $n_2=(\ell -1)/2$.

Using the lemma it follows that
$$
\frac{1}{E(q)^{\ell -6}} =\frac{(E(q)^3)^2}{E(q^\ell)} \equiv  \frac{(\sum_{i \in Jset(\ell)} J_i)^2} {E(q^\ell)} \pmod \ell
$$
and the powers of $q$  that are $r$ modulo $\ell$ do not show up. \halmos

It would be interesting to come up with an infinite family provable by Hirschhorn-style proofs!

{\bf Acknowledgment}: We are grateful to George Andrews, Bruce Berndt, Lev Borisov, Shaun Cooper, Frank Garvan,
and Michael Hirschhorn, for very useful advice, and to Silviu Radu for permission to post [R2].

{\bf References}

[B] Matthew Boylan, {\it Exceptional congruences for powers of the partition functions},
Acta Arith. {\bf 111} (2004), 187-203.

[D] Freeman J. Dyson, {\it Missed opportunities}, Bulletin of the American Mathematical Society {\bf 78} (1972), 635–652.

[Ha] G.H. Hardy, {\it ``Ramanujan''}, Cambridge University Press, 1940.

[Hi] Michael D. Hirschhorn, {\it A short and simple proof of Ramanujan's mod 11 partition congruence}, 
preprint available from \hfill\break
{\tt http://web.maths.unsw.edu.au/\~{}mikeh/webpapers/paper188.pdf}

[M] I.G. Macdonald, I. G., {\it Affine root systems and Dedekind's η-function}, Inventiones Mathematicae {\bf 15} (1972), 91–143

[R1] Silviu Radu, {\it An algorithmic approach to Ramanujan's congruences},  Ramanujan J. {\bf 20} (2009), 215-251 .

[R2] Silviu Radu, {\it Email message to Doron Zeilberger},  \hfill\break
 {\tt http://www.math.rutgers.edu/\~{}zeilberg/mamarim/mamarimhtml/SilviuRaduMessageJune2013.pdf}

[W] Lasse Winquist, {\it An elementary proof of} $p(11m+6) \equiv 0 \pmod {11}$, J. Combinatorial Theory {\bf 6}, 56–59.

\bigskip
\bigskip
\hrule
\bigskip

Edinah K. Gnang, Computer Science Department, Rutgers University (New Brunswick), Piscataway, NJ 08854, USA.
{\tt gnang at cs dot rutgers dot edu}

Doron Zeilberger, Mathematics Department, Rutgers University (New Brunswick), Piscataway, NJ 08854, USA.
{\tt zeilberg at math dot rutgers dot edu}

\end